\documentclass[a4paper,12pt,oneside]{article} 
\usepackage{amsfonts,amssymb}
\usepackage{color} \usepackage{mathrsfs} \let\mathcal\mathscr

\usepackage[all,ps,cmtip]{xy}

\title{\sc Singularities of the Projective Dual Variety}
\author{\sc Roland Abuaf \footnote{Institut Fourier, 100 rue des maths. Saint Martin d'Hères, 38402. E-mail :\it{abuaf@ujf-grenoble.fr}}}

\usepackage[dvips]{graphicx}

\usepackage{pstricks}
\usepackage{amsfonts,amssymb}
\usepackage{color} \usepackage{mathrsfs} \let\mathcal\mathscr

\usepackage[all,ps,cmtip]{xy}

\DeclareGraphicsExtensions{eps}
\DeclareGraphicsExtensions{pdf}

\usepackage[T1]{fontenc}
\usepackage{amssymb}

\usepackage{amsmath}
\usepackage{latexsym}

\usepackage[latin1]{inputenc}

\newtheorem{theo}{Theorem}[subsection]

\newtheorem{exem}[theo]{Example}
\newtheorem{rem}[theo]{Remark}
\newtheorem{prop}[theo]{Proposition}

\newtheorem{defi}[theo]{Definition}
\newtheorem{nota}[theo]{Notations}
\newtheorem{lem}[theo]{Lemma}
\newtheorem{cor}[theo]{Corollary}
\newtheorem{conj}[theo]{Conjecture}

\newenvironment{proof}
{
\noindent
\textit{\underline{Proof}} :\\
$\blacktriangleright\;$%
}
{\hspace{\stretch{1}}%
$\blacktriangleleft$}

\newenvironment{proofmain}
{
\noindent
\textit{\underline{Proof of the main theorem}} :\\
$\blacktriangleright\;$%
}
{\hspace{\stretch{1}}%
$\blacktriangleleft$}

{
\noindent
\textit{\underline{Road Map}} :\\
$\blacktriangleright\;$%
}
{\hspace{\stretch{1}}%
$\blacktriangleleft$}

{
\noindent
\textit{\underline{Sketch of the Proof}} :\\
$\blacktriangleright\;$%
}
{\hspace{\stretch{1}}%
$\blacktriangleleft$}

{
\noindent
\textit{\underline{Preuve}} : (id\'ee)\\
$\blacktriangleright\;$%
}
{\hspace{\stretch{1}}%
$\blacktriangleleft$}

\begin{document}

\maketitle

\begin{abstract}

Let $X \subset \mathbb{P}^N$ be an irreducible, non degenerate projective variety and let $X^* \subset {\mathbb{P}^N}^*$ be its projective dual. Let $L \subset \mathbb{P}^N$ be a linear space such that $\langle L,T_{X,x} \rangle \neq \mathbb{P}^N$ for all $x \in X_{smooth}$ and such that the lines in $X$ meeting $L$ do not cover $X$. If $x \in X$ is general, we prove that the multiplicity of $X^*$ at a general point of $\langle L,T_{X,x}\rangle ^{\bot}$ is strictly greater than the multiplicity of $X^*$ at a general point of $L^{\bot}$. This is a strong refinement of Bertini's theorem.

\end{abstract}

\vspace{\stretch{1}}

\newpage

\begin{section}{Introduction}

\subsection{Multiplicities of the Projective Dual}

\bigskip
Let $X \subset \mathbb{P}^N$ be an irreducible projective variety over the field of complex numbers. Let $X^* \subset {\mathbb{P}^N}^*$ be its projective dual, let $L \subset \mathbb{P}^N$ be a linear space and $H$ be a general hyperplane containing $L$. Bertini's classical theorem asserts that the tangency locus of $H$ with $X$ is included in $X\cap L$. Very little is known about the hyperplanes whose tangency locus with $X$ lies outside $L \cap X$. It is tempting to think that the multiplicity in $X^*$ of such a hyperplane is strictly larger than the multiplicity of a general hyperplane containing $L$. The following example shows that this is not true for every $L$.

\begin{exem} \label{scrollcubic} \upshape{Let $X \subset \mathbb{P}^4$ be a smooth hyperplane section of $\mathbb{P}^1 \times \mathbb{P}^2 \subset \mathbb{P}^5$. The variety $X$ is a ruled surface of degree $3$. Its dual is a hypersurface of degree $3$ in ${\mathbb{P}^4}^*$ which does not contain any points of multiplicity higher than $2$. Let $L$ be the exceptional section of $X$. If $H \subset \mathbb{P}^4$ is a general hyperplane which contains $L$, then $H \cap X = L \cup D_1 \cup D_2$, where $D_1$ and $D_2$ are two distinct lines on $X$ such that $D_1.D_2 = 0$ and $L.D_i =1$ for $i=1,2$. As a consequence, a general point of $L^{\bot}$ is of multiplicity $2$ in $X^*$. Now, let $D \subset X$ be a line such that $D.L = 1$ and let $x \in D$ such that $x \not\in L$. The hyperplane containing $L$ and $T_{X,x}$ is a point of multiplicity exactly $2$ in $X^*$, that is, the multiplicity of a general point of $L^{\bot}$.}
\end{exem}

This example shows that, even for general $x \in X$, the multiplicity in $X^*$ of a hyperplane containing $L$ and tangent to $X$ at $x$ may well be equal to the multiplicity of a general hyperplane containing $L$. Thus, without extra hypotheses on $L$, it seems hopeless to say something about the multiplicity in $X^*$ of special points of $L^{\bot}$. For this purpose, we introduce the following definition.

\begin{defi} \label{defidual} Let $X \subset \mathbb{P}^N$ be an irreducible projective variety and let $L \subset \mathbb{P}^N$ be a linear space. Consider the conormal diagram.

\begin{figure}[!h]
\centering
\label{conormal}
\begin{tabular}{rrcll}

& & $I(X/ \mathbb{P}^N) := \overline{ \{ (H,x) \in {\mathbb{P}^N}^* \times X_{smooth}, \,\, T_{X,x} \subset H \}} \subset {\mathbb{P}^N}^* \times \mathbb{P}^N$ & & \\
&  $\stackrel{q}\swarrow$ & & $\stackrel{p}\searrow$  &\\
&$X^* \subset {\mathbb{P}^N}^*$  & &  $ X \subset \mathbb{P}^N$&\\
\end{tabular}
\caption{conormal diagram}
\end{figure}

Let $F_1,...,F_m$ be all the irreducible components of $q^{-1}(L^{\bot})$ such that the restrictions:
$$ q|_{F_i} : F_i \rightarrow L^{\bot}$$ are surjective. The \emph{\textbf{contact locus}} of $L$ with $X$, which we denote by $\operatorname{Tan}(L,X)$, is the union of the $p(F_i)$, for $1 \leq i \leq m$.
\end{defi}

In the case where $L$ is a hyperplane, the contact locus $\operatorname{Tan}(L,X)$ is called the \textbf{\emph{tangency locus}} of $L$ with $X$. A \emph{\textbf{tangent hyperplane}} to $X$ is a hyperplane $H \subset \mathbb{P}^N$ such that $\operatorname{Tan}(H,X) \neq \emptyset$.

The contact locus $\operatorname{Tan}(L,X)$ can be thought as the variety covered by the tangency loci of general hyperplanes containing $L$. In case $L^{\bot} \not\subset X^*$, this locus is empty. We always have the inclusion:
$$ \overline{ \{ x \in X_{smooth}, \,\, T_{X,x} \subset L \}} \subset \operatorname{Tan}(L,X),$$
but if $\dim(L) < N-1$ or if $X$ is not smooth, the former locus can be strictly smaller than the latter. Note also that Bertini's theorem says that $\operatorname{Tan}(L,X) \subset L\cap X$. Finally, the contact locus is well behaved. If for a general hyperplane $H'$ containing $L$, we have $\dim \operatorname{Tan}(H',X) >0$, then:
$$ \operatorname{Tan}(H \cap L, H \cap X) = H \cap \operatorname{Tan}(L,X),$$ for any general hyperplane $H \subset \mathbb{P}^N$.

\begin{exem}
\upshape{If $X \subset \mathbb{P}^N$ is such that $X^*$ is a hypersurface and $L=T_{X,x}$, where $x \in X$ is a general point, then $\operatorname{Tan}(L,X) = x$. 

\noindent If $X = G(1,7) \subset \mathbb{P}^{27}$ and $L = \langle T_{X,y_1}, T_{X,y_2} \rangle$, where $y_1,y_2 \in \mathbb{G}(1,7)$ are two general points, then $\operatorname{Tan}(L,X) = \{ x \in X, \,\, T_{X,x} \subset L \}$ is a $4$-dimensional quadric, the entry locus of a general point $z \in \langle y_1,y_2 \rangle$.

\noindent If $X = G(1,4) \subset \mathbb{P}^9$ and $L = T_{X,y}$, for any $y \in X$, then $\dim \operatorname{Tan}(L,X) > 0$, whereas $\{ x \in X, \,\, T_{X,x} \subset L \} = \{ y \}$.
}
\end{exem}
\begin{defi}  Let $X \subset \mathbb{P}^N$ be an irreducible projective variety, and let $L \subset \mathbb{P}^N$ be a linear subspace. The \textbf{\emph{shadow}} of $L$ on $X$, which we denote by $\operatorname{Sh}_X(L)$, is the closed variety covered by the linear spaces $M \subset X$ such that $\dim(M) = \operatorname{def}(X)+1$ and $ \dim(M \cap \operatorname{Tan}(L,X)) = \operatorname{def}(X)$.
\end{defi}

Here $\operatorname{def}(X) = \operatorname{codim}(X^*)-1$. The shadow is also well behaved. Namely, assume that $\operatorname{def}(X) >0$, then:
$$ \operatorname{Sh}_L(X) = X \Leftrightarrow \operatorname{Sh}_{H\cap L}(H \cap X) = H \cap X,$$ for any general hyperplane $H \subset \mathbb{P}^N$.
Note also that if $x \in X$ is a general point and $L = T_{X,x}$, then $\operatorname{Sh}_L(X) \neq X$, unless $X$ is a linear space. Indeed, if $X^*$ is a hypersurface, this is obvious since $\operatorname{Tan}(T_{X,x},X)= x$ for general $x \in X$. If $X^*$ is not a hypersurface, take enough general hyperplane sections of $X$ passing through $x$, so that the corresponding dual is a hypersurface.

\noindent
Now we can state the main theorem of this paper.

\begin{theo} \label{maintheo} Let $X \subset \mathbb{P}^N$ be an irreducible, non-degenerate projective variety. Let $L \subset \mathbb{P}^N$ be a linear space such that $\operatorname{Sh}_X(L) \neq X$. Then, for all $x \in X_{smooth}$ such that $x \notin \operatorname{Sh}_X(L)$ and such that $ \langle L,T_{X,x}\rangle \neq \mathbb{P}^N$, the multiplicity in $X^*$ of a general hyperplane containing ${\langle L,T_{X,x}\rangle }$ is strictly larger than the multiplicity in $X^*$ of a general hyperplane containing $L$.
\end{theo}

If $X$ is the ruled cubic surface considered in example \ref{scrollcubic} and $L$ is the directrix of $X$, one notices easily that $\operatorname{Sh}_X(L)=X$. This shows that the hypothesis $\operatorname{Sh}_X(L) \neq X$ can not be withdrawn. An obvious corollary of theorem \ref{maintheo} is the following.

\begin{cor} Let $X \subset \mathbb{P}^N$ be an irreducible, non-degenerate projective variety. Let $L \subset \mathbb{P}^N$ be a linear space such that there is no line in $X$ which meet $L$. Then, for all $x \in X_{smooth}$ such that $x \notin L$ and $\langle L,T_{X,x}\rangle \neq \mathbb{P}^N$, the multiplicity in $X^*$ of a general hyperplane containing ${\langle L,T_{X,x}\rangle }$ is strictly larger than the multiplicity in $X^*$ of a general hyperplane containing $L$.
\end{cor}

\subsection{Variety of Multisecant Spaces and Duals}

\bigskip

We recall the definition of multisecant spaces to a projective variety.

\begin{defi}
Let $X \subset \mathbb{P}^N$ be an irreducible projective variety. Let
$$ {S^{k}_X}^0 = \left\{ (x_0,..,x_k,u) \in X\times ...\times X \times \mathbb{P}^N, \dim \langle x_0,...,x_k\rangle  = k, u \in \langle x_0,...x_k\rangle  \right\},$$ and let $S^k_X$ be its Zariski closure in $X \times...\times X \times \mathbb{P}^N.$

Denote by $\phi$ the projection onto $\mathbb{P}^N$. The variety $S^k(X) = \phi(S_X^k)$ is the $k$-th secant variety to $X$.
\end{defi}

\begin{theo}[Terracini's Lemma] Let $X \subset \mathbb{P}^N$ be an irreducible projective variety, and let $(x_0,...,x_k) \in X \times... \times X$, be general points. If $u$ is general in $\langle x_0,...,x_k \rangle$, we have the equality:
$$\langle T_{X,x_0},...,T_{X,x_k} \rangle = T_{S^k(X),u}.$$
\end{theo}

We refer to \cite{zak} for a proof.

\begin{defi} Let $X \subset \mathbb{P}^N$ be an irreducible, non-degenerate projective variety, and let $k$ be an integer such that $S^k(X) \neq \mathbb{P}^N$. We say that $X$ is \textbf{\emph{dual $k$-defective}} if $\operatorname{def}(S^k(X)) > t(S^k(X))$, where $t(S^k(X))$ is the dimension of the general fiber of the Gauss map of $S^k(X)$.
\end{defi}

Note that when $X$ is smooth, then dual $0$-defectivity is the classical dual defectivity. I don't know if there exist smooth varieties which are dual $k$-defective for some $k \geq 1$, but which are not dual $0$-defective. I believe it would be interesting to find some examples of such varieties. 

Note also that the notion of dual $k$-defectivity seems to be related to that of $R_k$ regularity explored by Chiantini and Ciliberto in \cite{cilichi}.

\bigskip

A consequence of the main theorem \ref{maintheo} and Terracini's lemma is the following :

\begin{prop} \label{dualsecant} Let $X \subset \mathbb{P}^N$ be an irreducible, non degenerate, smooth, projective variety. Assume moreover that for all $k$ such that $S^k(X) \neq \mathbb{P}^N$ the variety $X$ is not dual $k-1$-defective. Then, for any such $k$, we have:
$$ {S^k(X)}^* \subset X_{k+1}^*,$$
where $X_{k+1}^*$ is the set of points which have multiplicity at least $k+1$ in $X^*$.
\end{prop}

\begin{proof}
The case $k=0$ is the definition of ${S^0(X)}^* = X^*$. Let $k \geq 1$ be an integer such that $S^k(X) \neq \mathbb{P}^N$, let $z \in S^{k-1}(X)$ be a general point and $H$ be a general hyperplane containing $T_{S^{k-1}(X),z}$. Let's prove that:
$$\operatorname{Tan}(H,X) =  \{ x \in X, T_{X,x} \subset T_{S^{k-1}(X),z} \} .$$

Let $x_0,...,x_{k-1}$ be $k$ general points in $\operatorname{Tan}(H,X)$. Let $z'$ be a general point in $\langle x_0,...,x_{k-1} \rangle $, by Terracini's lemma we have:
$$T_{S^{k-1}(X),z'} = \langle T_{X,x_0},...,T_{X,x_{k-1}} \rangle.$$

So $z' \in \operatorname{Tan}(H,S^{k-1}(X))$. But by hypothesis, we have $\operatorname{def}(S^{k-1}(X)) = t(S^{k-1}(X))$, which implies that $$z' \in \overline{ \{ y \in S^{k-1}(X)_{smooth}, T_{S^{k-1}(X),y} = T_{S^{k-1}(X),z} \} },$$ so that $x_0,...,x_{k-1} \in  \{ x \in X, T_{X,x} \subset T_{S^{k-1}(X),z} \} $.

\bigskip
We now prove that $\operatorname{Sh}_X(T_{S^{k-1}(X),z}) \neq X$. The argument above shows that $$\operatorname{Tan}(T_{S^{k-1}(X),z},X) = \{x \in X, T_{X,x} \subset T_{S^{k-1}(X),z} \}.$$ Assume that $\operatorname{Sh}_X(T_{S^{k-1}(X),z}) = X$. Then, for all $x'' \in X$, there exists $x' \in \{x \in X, T_{X,x} \subset T_{S^k(X),z} \}$ such that the line $\langle x'',x' \rangle $ lies in $X$. But since $X$ is smooth, this line $ \langle x'',x' \rangle$ lies in $T_{X,x'}$. So we have $X \subset T_{S^{k-1}(X),z}$, which contradicts the non-degeneracy.

\bigskip

As a consequence of theorem \ref{maintheo}, we get that for a general $x \in X$, the multiplicity in $X^*$ of a general hyperplane containing $\langle T_{S^{k-1}(X),z},T_{X,x} \rangle$ is strictly larger than the multiplicity in $X^*$ of a general hyperplane containing $T_{S^{k-1}(X),z}$. We apply Terracini's lemma to find that $S^{k}(X)^* \subset X^*_{k+1}$. This concludes the proof.

\end{proof}
 
A stronger result than proposition \ref{dualsecant} has been stated for the first time by Zak in \cite{zak2}, but no proof was given there.

In the second part of this paper we present a proof of theorem \ref{maintheo}, while in the third part we discuss some consequences and open questions.

\bigskip

I would like to thank Bernard Teissier who took time to explain me part of his work and Christine Jost for interesting discussions on Segre classes. I am also especially grateful to Christian Peskine for the (numerous!) discussions we had together. This work owes a lot to his patience. Finally, my thanks go to the referee, who helped me to improve the paper, in many ways.

\newpage

\end{section}

\begin{section}{Proof of the Main Theorem}

When $Z \subset \mathbb{P}^N$, we denote by $\mathcal{C}_{z}(Z) \subset \mathbb{P}^N$ the embedded tangent cone to $Z$ at $z$ and if $H \subset \mathbb{P}^N$ is a hyperplane, then $[h]$ is the corresponding point in $(\mathbb{P}^N)^*$.

\bigskip

\noindent The proof of theorem \ref{maintheo} is obvious if $L^{\bot} \not\subset X^*$. Thus, we only deal with the case where $L^{\bot} \subset X^*$. Moreover, we can restrict to the case where $X^*$ is a hypersurface. Indeed, assume that $X^*$ has codimension $p\geq 2$. Let $z \in L^{\bot}$ and $z_x \in \langle L, T_{X,x} \rangle^{\bot}$ be general points, let $M \subset \mathbb{P}^N$ be a general $\mathbb{P}^{N+1-p}$ passing through $x$, let $X' = M \cap X$ and $L' = M \cap L$. We have $\operatorname{Sh}_{X'}(L') \neq X'$ and $\langle T_{X',x}, L' \rangle \neq \mathbb{P}^{N+1-p}$. Moreover, we have:
$$ (X')^* = \pi_{M^{\bot}}(X^*),$$
where $\pi_{M^{\bot}}$ is the projection from $M^{\bot}$ in ${\mathbb{P}^N}^*$. Since $M$ is general, the map $\pi_{M^{\bot}}$ is locally an isomorphism around $z_x$. Hence we have:
$$ \operatorname{mult}_{z}(X^*) = \operatorname{mult}_{z_x}(X^*) \Leftrightarrow \operatorname{mult}_{\pi_{M^{\bot}}(z)}((X')^*) = \operatorname{mult}_{\pi_{M^{\bot}}(z_x)}((X')^*).$$ Finally, note that $\pi_{M^{\bot}}(z)$ is a general point of $(L')^{\bot}$ and that $\pi_{M^{\bot}}(z_x)$ is a general point of $\langle L', T_{X',x} \rangle^{\bot}$. As a consequence, it is sufficient to prove the theorem for $X'$, whose dual is a hypersurface.

\bigskip

Let's start with a plan of the proof. We assume that $X^*$ has constant multiplicity along a smooth curve $S \subset L^{\bot}$ passing through ${\langle L, T_{X,x}\rangle }^{\bot}$ and through a general point of $L^{\bot}$ and we find a contradiction. More precisely:

\- We prove that the equimultiplicity of $X^*$ along $S$ implies that the family of the tangent cones to $X^*$ at the points of $S$ is flat.

\- Then, we show that the flatness of the family of the tangent cones to $X^*$ at the points of $S$ leads to the flatness of the family of the conormal spaces of these tangent cones. As a consequence, we have ${|\mathcal{C}_s(X^*)|}^* \subset L$ for all $s \in S$.

\- Finally, we relate the tangent cone to $X^*$ at $z$ to the set of tangent hyperplanes to $X^*$ at $z$ (when $z$ is a smooth point of $X^*$, this is the reflexivity theorem \cite{klei}). Using the fact that $\operatorname{Sh}_L(X) \neq X$, we deduce that $|\mathcal{C}_s(X^*)|^* \not\subset L$ for $s \in {\langle L, T_{X,x} \rangle}^{\bot}$ and thus a contradiction.

\begin{subsection}{Normal Flatness and Lagrangian Specialization Principle}

\bigskip

Let $S \subset Z \subset \mathbb{P}^N$ be two varieties. We recall some properties of the tangent cones $\mathcal{C}_{s}(Z), s \in S$ when $Z$ is equimultiple along $S$.

\begin{defi} Let $S \subset Z$ be two varieties. We say that $Z$ is \textbf{\emph{equimultiple}} along $S$ if the multiplicity of the local ring $\mathcal{O}_{Z,s}$ is constant for $s \in S$.
\end{defi}

\begin{prop}[\cite{hiro}, cor. 2, p. 197] Let $Z \subset \mathbb{P}^N$ be a hypersurface and $S$ a connected smooth subvariety (not necessarily closed) of $Z$ such that $Z$ is equimultiple along $S$.

Then, for all $s \in S$, there exists an open neighborhood $U$ of $s$ in $S$ containing $s$ and a closed subscheme $\mathcal{G}(Z) \subset \mathbb{P}^N \times U$ such that the  natural projection $p :\mathcal{G}(Z)~\rightarrow~U$ is a flat and surjective morphism whose fiber $\mathcal{G}(Z)_{s'}$ over any $s' \in U$ is $\mathcal{C}_{s'}(Z)$.
\end{prop}

We assume that our theorem is not true, that is for general $x \in X$, the multiplicity of $X^*$ at a general point of $\langle L,T_{X,x}\rangle ^{\bot}$ is equal to the multiplicity at a general point of $L^{\bot}$.

Let $[h]$ be a general point of $\langle L,T_{X,x}\rangle ^{\bot}$ and let $S \subset L^{\bot}$ be a smooth (not necessarily closed) connected curve passing through $[h]$ and through a general point of $L^{\bot}$. We apply the above proposition to $X^*$ and $S$. Then there exists a scheme $\mathcal{G}(X^*) \subset {\mathbb{P}^N}^* \times S$ such that the natural projection $p : \mathcal{G}(X^*) \rightarrow S$ is a flat and surjective morphism whose fiber over $s \in S$ is the tangent cone to $X^*$ at $s$. Let $\Gamma(X^*) = |\mathcal{G}(X^*)|$. The induced morphism $\Gamma(X^*) \rightarrow S$ is flat and for general $s \in S$ the fiber $\Gamma(X^*)_s$ is exactly $|\mathcal{C}_s(X^*)|$.
\bigskip

Now we study the family of the duals of the reduced tangent cones of $X^*$ at points of $S$. Applying the Lagrangian specialization principle (see \cite{lete} and \cite{klei2}) to $\Gamma(X^*)$ and $S$, we find the following.
\begin{theo} Let $S \subset X^*$ be a smooth curve such that $X^*$ is equimultiple along $S$. There esists a variety $I_S(\Gamma(X^*) / {\mathbb{P}^N}^* \times S)$ with the following properties.

i) For general $s \in S$, the following equality holds in $\mathbb{P}^N \times \Gamma(X^*)_s$:
$$I(|\mathcal{C}_s(X^*)| / {\mathbb{P}^N}^*) = I_S(\Gamma(X^*) / {\mathbb{P}^N}^* \times S)_s.$$

ii)The morphism $I_S(\Gamma(X^*) / {\mathbb{P}^N}^* \times S) \rightarrow S$ is flat and surjective,

iii)For all $s \in S$, the conormal space $I(|\mathcal{C}_s(X^*)|/ {\mathbb{P}^N}^*)$ is a union of irreducible components of the reduced fiber $|I_S(\Gamma(X^*) / {\mathbb{P}^N}^* \times S)_s|$.
\end{theo}

 As a consequence of the above theorem, the image in $\mathbb{P}^N$ of the fiber $I_{S}(\Gamma(X^*)/ {\mathbb{P}^N}^* \times S)_{s}$ is $|\mathcal{C}_s(X^*)|^*$, for general $s \in S$. Moreover, for any $s \in S$, the image of the reduced fiber $|I_S(\Gamma(X^*) / {\mathbb{P}^N}^*\times S)_s|$ contains $|\mathcal{C}_s(X^*)|^*$.

\end{subsection}

\begin{subsection}{Polar Varieties and Duals of Tangent Cones}

We discuss an extension of the reflexivity theorem proved by L\^e and Teissier in \cite{lete}. The main result of this section will be applied to $X^*$, so that we restrict our study to the case of hypersurfaces.

\begin{defi} Let $Z \subset \mathbb{P}^N$ be a reduced and irreducible hypersurface and let $D \subset \mathbb{P}^N$ be a linear space. The \textbf{\emph{polar variety }}of $Z$ associated to $D$, which we denote by $P(Z,D)$, is the closure of the set $\{z \in Z_{smooth}, D \subset T_{Z,z} \}$.

If $D = \emptyset$ (that is $D$ has dimension $-1$), then we put $P(Z,D) = Z$.
\end{defi}

\begin{rem} \upshape{ Note that if $Z$ is normal, if $u = [u_0,...,u_N]$ in an homogeneous system of coordinates on $\mathbb{P}^N$ and $f$ is an equation of $Z$ in this system then $P(Z,u)$ is given by the equations $f=0$ and $u_0 \frac{\partial f}{\partial x_0}+...+ u_N \frac{\partial f}{\partial x_N} = 0$.

If $Z$ is not normal, then all irreducible components of $Z_{sing}$ which are of dimension $N-2$ are irreducible components of the scheme defined by $f=0$ and $u_0 \frac{\partial f}{\partial x_0}+...+ u_N \frac{\partial f}{\partial x_N} = 0$, but they are not irreducible components of $P(X,u)$.}

\end{rem}

\begin{prop} Let $Z \subset \mathbb{P}^N$ be a reduced, irreducible hypersurface and let $D \subset \mathbb{P}^N$ be a general linear space of dimension $k$. Then $P(Z,D)$ is empty or of codimension $k+1$ in $Z$.
\end{prop}

\bigskip

We state a result of L\^e and Teissier which relates the duals of the tangent cones at $z$ of some polar varieties of $Z$ with the tangency locus of $z^{\bot}$ with $Z^*$. See $\cite{lete}$, proposition 2.2.1. 
For any $z \in Z$, recall that $\operatorname{Tan}(z^{\bot},Z^*)$ is the tangency locus of $z$ along $Z^*$ (see definition \ref{conormal}).

\begin{theo} \label{lete} Let $Z \subset \mathbb{P}^N$ be a reduced and irreducible hypersurface and let $z \in Z$ be a point. We have the following:

i) The dual of $|\mathcal{C}_z(Z)|$ is a union of reduced spaces underlying (possibly embedded) components of $\operatorname{Tan}(z^{\bot},Z^*)$.

ii) Any irreducible component of $|\operatorname{Tan}(z^{\bot},Z^*)|$ is dual to an irreducible component of $|\mathcal{C}_z(P(Z,D))|$ for general $D \in \mathbb{G}(k,N)$ and for some integer $k \in \{-1,...,N-2 \}$.

\end{theo}

\begin{rem} \upshape{ The ii) of the theorem has to be explained. Assume that there is an irreducible component (say $T$) of $|\operatorname{Tan}(z^{\bot},Z^*)|$ which is not dual to an irreducible component of $|\mathcal{C}_z(Z)|$. Then, there is $k \in \{0,...,N-2 \}$ such that for general $D \in \mathbb{G}(k,N)$, we have $z \in P(Z,D)$. Moreover, as $D$ varies in a dense open subset of $\mathbb{G}(k,D)$, the cones $\mathcal{C}_z(P(D,Z))$ have a fixed irreducible component in common whose reduced locus is $T^*$.

Note also that if $z \in Z_{smooth}$ then for $k \geq 0$ and for $D$ general in $\mathbb{G}(k,N)$, we have $z \notin P(Z,D)$. As a consequence of the ii) of the above theorem, we find $\operatorname{Tan}(z^{\bot},Z^*) = T_{Z,z}^{\bot}$ for $z \in Z_{smooth}$. This is the way the (obvious corollary of the) reflexivity theorem is often stated.}
\end{rem}

When $\operatorname{Tan}(z^{\bot},Z^*)$ is irreducible, one may expect $|\mathcal{C}_z(Z)|^* = |\operatorname{Tan}(z^{\bot},Z^*)|$. The following example shows that it is not true.

\begin{exem} \upshape{Let $X \subset \mathbb{P}^4$ be the smooth ruled surface of degree $3$ considered in example 1.1.1 and let $X^*$ its dual. The hypersurface $X^*$ has also degree $3$ and its singular locus is a $\mathbb{P}^2$, the dual of the exceptional section of $X$ (which we denote by $L$). Let $C \subset L^{\bot} = X^*_{sing}$ be the conic corresponding to the hyperplanes which are tangent to $X$ along a ruling of $X$ and let $z \in C$.

The tangent cone $\mathcal{C}_{z}(X^*)$ is a doubled $\mathbb{P}^3$ so that $|\mathcal{C}_z(X^*)|^* \neq \operatorname{Tan}(z^{\bot},X)$. We also note that the scheme-theoretic tangency locus of $z^{\bot}$ along $X$ is a line with an embedded point. The embedded point is dual to $|\mathcal{C}_z(X^*)|$ and the line is dual to $|\mathcal{C}_z(P(X^*,u))|$, for general $u \in {\mathbb{P}^4}^*$.}
\end{exem}

Before giving the proof of theorem \ref{maintheo}, we need some handy notations.

\begin{nota} Let $f : Y \rightarrow T$ be a quasi-projective morphism between quasi-projective schemes, let $T' \subset T$ be a smooth variety and let $s \in T'$ be any point. Let $Y_1,...,Y_m$ be the irreducible components of $f^{-1}(T')$ such that the restrictions:
$$ f|_{Y_i} : Y_i \rightarrow T',$$ are surjective. We denote by $\operatorname{limflat}_{ \{ t \rightarrow s, t \in T' \} } f^{-1}(t)$ the scheme:
$$\operatorname{limflat}_{\{t \rightarrow s, t \in T' \} } f^{-1}(t) = f|_{Y_1 \cup ... \cup Y_m}^{-1}(s).$$ If $\dim (T') = 1$ and the $Y_i$ are all reduced then $\operatorname{limflat}_{ \{t \rightarrow s, t \in T' \} } f^{-1}(t)$ is the classical flat limit taken along a smooth curve. If $f|_{f^{-1}(T')} : f^{-1}(T') \rightarrow T'$ is flat, then:
$$\operatorname{limflat}_{\{t \rightarrow s, t \in T' \} } f^{-1}(t) = f|_{f^{-1}(T')}^{-1}(s).$$
\end{nota}
\bigskip

\begin{proofmain}
We recall the setting for the convenience of the reader. The projective variety $X \subset \mathbb{P}^N$ is irreducible and non degenerate. The linear space $L \subset \mathbb{P}^N$ is such that $\operatorname{Sh}_X(L) \neq X$ and $\langle L,T_{X,x} \rangle \neq \mathbb{P}^N$ for all $x \in X_{smooth}$. We want to prove that for all $x \in X_{smooth}$ such that $x \notin \operatorname{Sh}_X(L)$, the multiplicity in $X^*$ of a general hyperplane containing $\langle L,T_{X,x} \rangle$ is strictly greater than that of a general hyperplane containing $L$. 
\bigskip

The result is obvious if $L^{\bot} \not\subset X^*$ and we have already seen that we can restrict to the case where $X^*$ is a hypersurface. So we only consider the case where $L^{\bot} \subset X^*$ and $X^*$ is a hypersurface and we assume that our result is not true. Let $x \in X_{smooth}$ with $x \notin \operatorname{Sh}_X(L)$ and let $[h]$ be a general point in $\langle L,T_{X,x} \rangle^{\bot}$. By the results of the previous section, there exists a smooth (non necessarily closed) curve $S \subset L^{\bot}$ with $[h] \in S$ and a flat morphism:
$$I_S(\Gamma(X^*) / {\mathbb{P}^N}^* \times S) \rightarrow S,$$ whose fiber $I_S(\Gamma(X^*) / {\mathbb{P}^N}^* \times S)_s$ is the the conormal space of $|\mathcal{C}_s(X^*)|$, for general $s \in S$. Moreover, the conormal space of $|\mathcal{C}_s(X^*)|$ is included in $|I_S(\Gamma(X^*) / {\mathbb{P}^N}^* \times S)_s|$ for all $s \in S$.

\bigskip
 
\noindent The i) of theorem \ref{lete} implies:
$$ |\mathcal{C}_s(X^*)|^* \subset p(|q^{-1}(s)|),$$ for all $s \in S$, where $p$ and $q$ are defined in the conormal diagram (see \emph{Figure 1}). Moreover, the flatness of $I_S(\Gamma(X^*) / {\mathbb{P}^N}^* \times S) \rightarrow S$ gives the inclusion:
$$ |\mathcal{C}_{[h]}(X^*)|^* \subset p(\operatorname{limflat}_{ \{ s \rightarrow [h], s \in S \} }|q^{-1}(s)|). $$ By definition \ref{defidual}, we have $p(\operatorname{limflat}_{ \{ s \rightarrow [h], s \in S \} }|q^{-1}(s)|) \subset \operatorname{Tan}(L,X) \subset L$.

\bigskip

\noindent Let $\mathcal{F}$ be an irreducible component of $\operatorname{Tan}(H,X)$ passing through $x$. By Theorem \ref{lete}, there is an integer $k \in \{-1,...,N-2\}$ such that $|\mathcal{F}|$ is dual to an irreducible component of $|\mathcal{C}_{[h]}(P(X^*,D))|$, for general $D \in \mathbb{G}(k,N)$. Since $|\mathcal{C}_{[h]}(X^*)|^* \subset L$, we have $k \geq 0$.

\bigskip

\noindent Let $x_0 \in \mathcal{F}$ be a general point. Duality implies $T_{|\mathcal{C}_{[h]}(P(X^*,D))|,z} \subset x_0^{\bot}$ for some general $z$ in the irreducible component of $\mathcal{C}_{[h]}(P(X^*,D))$ whose reduced locus is $|\mathcal{F}|^*$. Note that $\mathcal{C}_{[h]}(P(X^*,D)) \subset \mathcal{C}_{[h]}(X^*)$. Let $T_{|\mathcal{C}_{[h]}(X^*)|,z}$ be a limit of tangent spaces to $|\mathcal{C}_{[h]}(X^*)|$ at $z$. The point $z$ is general in $|\mathcal{C}_{[h]}(P(X^*,D))|$, so $ T_{|\mathcal{C}_{[h]}(P(X^*,D))|,z} \subset T_{|\mathcal{C}_{[h]}(X^*)|,z}$. 

\bigskip

\noindent As a consequence of this, we have $T_{|\mathcal{C}_{[h]}(P(X^*,D))|,z} \subset x_0^{\bot} \cap T_{\mathcal{C}_{[h]}(X^*),z}$. That is:
$$\langle x_0, T_{|\mathcal{C}_{[h]}(X^*)|,z}^{\bot} \rangle \subset \mathcal{F} \subset X.$$ But we have $|\mathcal{C}_{[h]}(X^*)|^* \subset \operatorname{Tan}(L,X)$, so that $T_{|\mathcal{C}_{[h]}(X^*)|,z}^{\bot} \in \operatorname{Tan}(L,X)$ and the above inclusion says that $x_0 \in \operatorname{Sh}_X(L)$. This is a contradiction.

\end{proofmain}

\end{subsection}

\end{section}

\newpage

\begin{section}{Corollaries and Open Questions}
We present here some corollaries of theorem \ref{maintheo} and related open questions.

\begin{subsection}{Zak's Conjecture on Varieties with Minimal Codegree}

\bigskip

Let $X \subset \mathbb{P}^N$ be an irreducible, non-degenerate projective variety. We recall, following Zak, that the order of $X$ is $\operatorname{ord}{X} = \min \{k, S^{k-1}(X) = \mathbb{P}^N \}$ and the $k$-th secant-defect is $\delta_k = \dim X + \dim S^{k-1}(X) +1 - \dim S^k(X)$, for all $k\leq \operatorname{ord}(X)-1$.

In \cite{zak}, Zak proves an important result related to secant defects.

\begin{theo}[Zak's Superadditivity Theorem]
Let $X \subset \mathbb{P}^N$ an irreducible, non-degenerate projective variety such that $\delta_1 >0$. For all $k \leq \operatorname{ord}(X)-1$, we have the inequality:
$$ \delta_k \geq \delta_{k-1} + \delta_1.$$
\end{theo}

The varieties on the boundary are called Scorza varieties, more precisely:

\begin{defi} An irreducible, smooth, non-degenerate projective variety $X \subset \mathbb{P}^N$ is a Scorza variety if the following conditions hold:

i) $\delta_1 >0$ and $N > 2n+1- \delta_1$,

ii) $\delta_k = \delta_{k-1} + \delta_{1}$ for all $k \leq \operatorname{ord}(X)-1$,

iii) $\operatorname{ord}(X)-1 =  [\frac{\dim X}{\delta_1}] $, where $[.]$ is the integral part.
\end{defi}

Zak gives in \cite{zak} a classification of Scorza varieties.

\begin{theo}[Classification of Scorza Varieties] Let $X$ a Scorza variety, then $X$ is one of the following:

\begin{tabular}{l p{11cm}}
 (i) & $X = v_2(\mathbb{P}^n) \subset \mathbb{P}^{n(n+3)/2}$ ($2^{nd}$ Veronese) and $\deg(X^*)=n+1$;\\
 (ii) & $X=\mathbb{P}^{n}\times\mathbb{P}^{n} \subset \mathbb{P}^{n(n+2)}$ and $\deg(X^*)=n+1$;\\
 (iii) & $X = \mathbb{G}(1,2n+1) \subset \mathbb{P}(\Lambda^2 \mathbb{C}^{2n+2})$ and $deg(X^*)=n+1$;\\
 (iv) & $X \subset \mathbb{P}^{26}$ is the $16$-dimensional variety corresponding to the orbit of highest weight vector in the lowest non trivial representation of the group of type $E_6$ and $\deg(X^*)=3$.\\
 \end{tabular}

\end{theo}

Zak notices in \cite{zak2} an important consequence of the assertion $S^k(X)^* \subset X^*_{k+1}$ (where $X^*_k$ is the set of points of multiplicity at least $k$ in $X^*$). We state his result in the setting where we are able to prove it.

\begin{prop} Let $X \subset \mathbb{P}^N$ be an irreducible, non-degenerate, smooth, projective variety. Assume that $X$ is not $k$ dual defective for $k < \operatorname{ord}(X)-1$, then the following inequality holds:
$$ \deg(X^*) \geq \operatorname{ord}(X).$$
\end{prop}

\begin{proof} With the above assumptions, proposition 1.2.4 implies that there is a point of multiplicity $\operatorname{ord}(X)-1$ in $X^*$. Since $X$ is non degenerate, its dual is not a cone and so $\deg(X^*) \geq \operatorname{ord}(X)$.

\end{proof}

\bigskip

If $X$ is a Scorza variety then $\deg(X^*) = \operatorname{ord}(X)$. Zak conjectures in \cite{zak2} the converse statement. We formulate his conjecture in the setting where we can prove the inequality: $\deg(X^*) \geq \operatorname{ord}(X)$.

\begin{conj}[\cite{zak2}] Let $X \subset \mathbb{P}^N$ be an irreducible, smooth, non-degenerate, projective variety. Assume that $X$ is not $k$ dual defective for all $k < \operatorname{ord}(X)$ and that $\deg(X^*) = \operatorname{ord}(X)+1$, then $X$ is a hyperquadric or a Scorza variety.

\end{conj}

It is proved in \cite{zak}, without any hypothesis on the dual defectiveness of $X$, that smooth varieties with $\operatorname{deg}(X^*) = 3$ and $\operatorname{ord}(X)=3$ are Severi varieties. In particular, they are Scorza varieties. Note, however, that the smoothness assumption seems to be necessary in his proof. I believe it would be very interesting to have a classification of all varieties whose duals have degree $3$.

\end{subsection}

\begin{subsection}{Varieties with Unexpected Equisingular Linear Spaces}

We come back to our usual setting. Let $L \subset \mathbb{P}^N$ be a linear space such that for all $x \in X_{smooth}$, we have $\langle L,T_{X,x} \rangle \neq \mathbb{P}^N$. We have seen in example 1.1.1 that a hyperplane containing the join $\langle L, T_{X,x} \rangle$ may have the same multiplicity in $X^*$ as the general hyperplane containing $L$, even if $x$ is a general point of $X$. The following definition is convenient to describe this situation.

\begin{defi} Let $X \subset \mathbb{P}^N$ be an irreducible, non-degenerate projective variety such that $X^*$ is a hypersurface. Let $L \subset \mathbb{P}^N$ be a linear space such that for all $x \in X_{smooth}$, we have $\langle L, T_{X,x} \rangle \neq \mathbb{P}^N$. We say that $L^{\bot}$ is an \textbf{\emph{unexpected equisingular linear space}} in $X^*$ if for all $x \in X_{smooth}$, the general hyperplane containing $\langle L,T_{X,x} \rangle$ has the same multiplicity in $X^*$ as the general hyperplane containing $L$.
\end{defi}

The variety in example \ref{scrollcubic} is rather special since it is a scroll surface (see \cite{zak2} for interesting discussions about this variety). It is not a coincidence that the directrix of this variety is an unexpected equisingular linear space in its dual. Indeed, we have the following result.

\begin{theo} \label{converse} Let $X \subset \mathbb{P}^N$ be an irreducible, smooth, non-degenerate projective variety such that $X^*$ is a hypersurface. Let $L \subset X$ be a linear space with $\dim(L) = \dim(X)-1$. Assume that $L^{\bot}$ is an unexpected equisingular linear space in $X^*$ such that $\operatorname{mult}_{L^{\bot}}(X^*) =2$. Then $X$ is the cubic scroll surface in $\mathbb{P}^4$.

\end{theo}

Here $\operatorname{mult}_{L^{\bot}}(X^*)$ denotes the multiplicity in $X^*$ of a general point of $L^{\bot}$. Before diving into the proof of theorem \ref{converse}, we describe the tangency locus of any point $[h] \subset X^*$, such that $\operatorname{mult}_{[h]}(X^*) = 2$.

\begin{prop} \label{tangency} Let $X \subset \mathbb{P}^N$ be a smooth, irreducible, non-degenerate projective variety such that $X^*$ is a hypersurface. Let $[h] \in X^*$ be such that $\operatorname{mult}_{[h]}(X^*) = 2$. The scheme theoretic tangency locus of $H$ with $X$ is either:

i) an irreducible hyperquadric and in this case $|\mathcal{C}_{[h]}(X^*)|^* = \operatorname{Tan}(H,X)$,

ii) the union of two (not necessarily distinct) linear spaces,

iii) a linear space with at least one embedded component.
\end{prop} 
 
We postpone the proof of this proposition to the appendix, and we start the proof of theorem \ref{converse}.

\begin{proof}

Let $H$ be a general hyperplane containing $L$. We have $H \cap X = L \cup D_H$, where $D_H$ is a divisor such that:
$$ D_H \cap L = \operatorname{Tan}(H,X).$$
Let $x \in X$ be a general point and let $H_x$ be a general hyperplane containing $\langle L,T_{X,x} \rangle$. Then $\operatorname{Tan}(H_x,X)$ contains $x$ and $p(\operatorname{limflat}_{ \{ [h] \rightarrow [h_x], [h] \in L^{\bot} \} }q^{-1}([h]))$.
By hypothesis, we have:
$$ \operatorname{mult}_{[h_x]}(X^*) = \operatorname{mult}_{[h]}(X^*) = 2,$$ for all $[h] \in L^{\bot}$. Proposition \ref{tangency} hence implies that the irreducible component of $\operatorname{Tan}(H_x,X)$ containing $x$, which we denote by $R_{H_x}$, also contains $p(\operatorname{limflat}_{ \{ [h] \rightarrow [h_x], [h] \in L^{\bot} \} }q^{-1}([h]))$. Moreover, we have:
$$p(\operatorname{limflat}_{ \{ [h] \rightarrow [h_x], [h] \in L^{\bot} \} }q^{-1}([h])) \subset L,$$ so that $$\dim R_{H_x} > \dim p(\operatorname{limflat}_{ \{ [h] \rightarrow [h_x], [h] \in L^{\bot} \} }q^{-1}([h])),$$ for general $[h] \in L^{\bot}$. As a consequence $\dim R_{H_x} = n-1$.
 
\bigskip

On the other hand, since 
$$ \operatorname{mult}_{[h_x]}(X^*) = \operatorname{mult}_{[h]}(X^*) = 2,$$ for all $[h] \in L^{\bot}$, we have $|\mathcal{C}_{[h_x]}(X^*)|^* \neq |R_{H_x}|$. We apply again proposition $\ref{tangency}$ and we find that $|R_{H_x}|$ is necessarily a linear space of dimension $n-1$. Thus, we have:
$$ \dim \langle L, T_{X,x} \rangle = n+1.$$Note that Bertini's theorem implies:
$$R_{H_x} \subset \langle L, T_{X,x} \rangle \cap X,$$ for general $H_x$ containing $\langle L, T_{X,x} \rangle$.  As a consequence $R_{H_x}$ is an irreducible component of $\langle L, T_{X,x} \rangle \cap X$, for general $H_x$. Thus $R_{H_x}$ does not depend of $H_x$, for general $H_x$ containing $\langle L, T_{X,x} \rangle$. We deduce that $\langle L, T_{X,x} \rangle $ is tangent to $X$ along a linear space of dimension $n-1$. By the theorem on tangencies, we have $n-1 \leq 1$, that is $n=2$ (obviously, $X$ is not a curve). So $X \subset \mathbb{P}^N$ is a non degenerate surface containing a distinguished line $L$, such that for general $x \in X$, there is a $\mathbb{P}^3$ tangent to $X$ along a line passing through $x$ and meeting $L$. This means that $X$ is the projection of a scroll of type $S_{1,d-1}$. By hypothesis, we have $\operatorname{mult}_{L^{\bot}}(X^*) = 2$, hence proposition 1.6. of \cite{cilirussosimis} implies that $X = S_{1,2} \subset \mathbb{P}^4$.

\end{proof}

\end{subsection}

\end{section}

\newpage

\begin{appendix}
\begin{section}{Tangency Loci of Points of Multiplicity $2$ in the Dual}
The goal of this appendix is to prove the following proposition.

\begin{prop} \label{multideux} Let $X \subset \mathbb{P}^N$ be a smooth, irreducible, non-degenerate projective variety such that $X^*$ is a hypersurface. Let $[h] \in X^*$ be such that $\operatorname{mult}_{[h]}(X^*) = 2$. The scheme theoretic tangency locus of $H$ with $X$ is either:

i) an irreducible hyperquadric and in this case $|\mathcal{C}_{[h]}(X^*)|^* = \operatorname{Tan}(H,X)$,

ii) the union of two (not necessarily distinct) linear spaces,

iii) a linear space with at least one embedded component.
\end{prop} 

\begin{exem} \upshape{All three cases can be encountered in Nature.

\noindent i) If $X = v_2(\mathbb{P}^2) \subset \mathbb{P}^5$, then for all $[h] \in v_2({\mathbb{P}^2}^*) \subset X^*$, we have $\operatorname{mult}_{[h]}(X^*) = 2$ and $\operatorname{Tan}(H,X)$ is a smooth conic.

\noindent ii) If $X$ is a complete intersection of large multidegree and large codimension, then there are points $[h_1], [h_2] \in X^*$ such that $\operatorname{mult}_{[h_i]}(X^*) = 2$ and $\operatorname{Tan}(H_1,X)$ is exactly two distinct points, whereas $\operatorname{Tan}(H_2,X)$ is a single double point.

\noindent iii) If $X$ is the cubic scroll of example \ref{scrollcubic}, then there is a conic $C \subset X^*$, such that for all $[h] \in C$, we have $\operatorname{mult}_{[h]}(X^*) = 2$ and $\operatorname{Tan}(H,X)$ is a line with an embedded point.

}

\end{exem}

A doubled linear space will be considered as the union of two (not distinct) linear spaces. By theorem \ref{lete}, we know that the irreducible components of $\operatorname{Tan}(H,X)$ are dual to irreducible components of the reduced spaces underlying some $\mathcal{C}_{[h]}(P(X^*,D_k))$ for general $D_k \in \mathbb{G}(k,N)$. In the case where $\operatorname{mult}_{[h]}(X^*) = 2$, the cones $\mathcal{C}_{[h]}(P(X^*,D_k))$ are rather easy to describe. Let's start with some notations.

\begin{nota}
Let $Z \subset \mathbb{P}^N$ be a reduced and irreducible hypersurface. Let $D \in \mathbb{G}(k,N)$ and let $f_Z$ be an equation for $Z$ in some coordinate system of $\mathbb{P}^N$. We denote by $P(f_Z,D)$ the subscheme of $\mathbb{P}^N$ whose ideal is generated by the equations: $$u_0\frac{\partial f_Z}{\partial t_0} +...+ u_N\frac{\partial f_Z}{\partial t_N},$$
for $u = [u_0,...,u_N]$ varying in $D$.
\end{nota}

Let $D \in \mathbb{G}(k,N)$ be a general $k$-plane. Note that if $\dim(Z_{sing}) < \dim P(Z,D)$ (that is $\dim Z_{sing} \leq N-k-3$), then $P(Z,D) = P(f_Z,D) \cap Z$. In the other case, the irreducible components of maximal dimension of $Z_{sing}$ are irreducible components of $P(f_Z,D) \cap Z$.

\begin{lem} Let $Z \subset \mathbb{P}^N$ be an irreducible and reduced hypersurface. Let $z \in Z$ and let $k \in\{-1,...,N-2\}$. Then, for general $D \in \mathbb{G}(k,N)$, we have:

1) $ z \notin P(Z,D)$ or,

2) $\operatorname{mult}_z P(Z,D) = \operatorname{mult}_z(Z).\operatorname{mult}_z P(f_Z,D)$, if  $\dim(Z_{sing}^{(z)}) < \dim P(Z,D)$, where $Z_{sing}^{(z)}$ is an irreducible component of $Z_{sing}$ of maximal dimension passing through $z$.

3) $\operatorname{mult}_{z} P(Z,D) < \operatorname{mult}_z(Z).\operatorname{mult}_z P(f_Z,D)$, if $\dim(Z_{sing}^{(z)}) \geq \dim P(Z,D)$,
where $Z_{sing}^{(z)}$ is an irreducible component of $Z_{sing}$ of maximal dimension passing through $z$.
\end{lem}

\begin{proof} If $z \in P(Z,D)$ for general $D \in \mathbb{G}(k,N)$, we will prove the lemma only in the case $P(f_Z,D)$ is smooth at $z$, for two reasons. The general case is obtained by the same methods, this is only more technical, and we will use the result only in the case $P(f_Z,D)$ is smooth at $z$.

Moreover if $z \in P(Z,D)$ for general $D$, we will only concentrate on the case $\dim (Z_{sing}^{(s)}) < \dim P(Z,D)$. In this case, we have locally around $z$ the equality $P(Z,D) = P(f_Z,D) \cap Z$ for genral $D \in \mathbb{G}(k,N)$. The situation where an irreducible component $Z_{sing}$ containing $z$ is an irreducible component of $P(f_Z,D) \cap Z$ (which is case 3 of the lemma) is dealt with exactly in the same way.

\bigskip
Now, we work locally around $z$, so that $P(f_Z,D) \cap Z = P(Z,D) \subset \mathbb{A}^N$, for general $D \in \mathbb{G}(k,N)$.
Let $(Z_i)_{i \in I}$ be a stratification of $Z$ such that $Z_i$ is smooth and $Z$ is normally flat along $Z_i$, for all $i \in I$. Such a stratification exists, due to the open nature of normal flatness (see \cite{hiro}, chapter II). Consider the Gauss map $G : Z \rightarrow {\mathbb{P}^N}^*$. It restricts to a map $G_i : Z_i \rightarrow {\mathbb{P}^N}^*$. We have: 
$$P(f_Z,D) \cap Z = P(Z,D) = G^{-1}(D^{\bot}),$$ so that
$P(f_Z,D) \cap Z_i = G_i^{-1}(D^{\bot})$, for all $i$.

Now, we apply Kleiman's transversality theorem to find that for all $i$ and for general $D \in \mathbb{G}(k,N)$, the inverse images $G_i^{-1}(D^{\bot})$ are either empty or smooth of the expected dimension.

Let $i$ such that $z$ is in $Z_i$. If $z \notin G_i^{-1}(D^{\bot})$ for general $D \in \mathbb{G}(k,N)$, then $z \notin P(Z,D)$ and we are in the case 1 of the lemma. Otherwise, $z$ is a smooth point of $G_i^{-1}(D^{\bot})$, so $T_{P(f_Z,D),z}$ and $T_{Z_i,z}$ are transverse.

\bigskip

Assume that $\operatorname{mult}_z P(Z,D) > \operatorname{mult}_z(Z).\operatorname{mult}_z P(f_Z,D)$. Since $P(f_Z,D)$ is smooth at $z$, this implies that $T_{P(f_Z,D),z}$ and $\mathcal{C}_z(Z)$ are not transverse. In particular, the linear spaces $T_{P(f_Z,D),z}$ and $\operatorname{Vert}(\mathcal{C}_z(Z))$ are not transverse (here $\operatorname{Vert}(\mathcal{C}_z(Z))$ is the vertex of the cone $\mathcal{C}_z(Z)$). But $Z$ is normally flat along $Z_i$, so we have $T_{Z_i,z} \subset \operatorname{Vert}(\mathcal{C}_z(Z))$ (see thm. 2, p. 195 of \cite{hiro}). This is a contradiction.
\end{proof}

As a consequence, we have the following corollary.

\begin{cor} \label{cone} Let $Z \subset \mathbb{P}^N$ be a reduced, irreducible hypersurface. Let $z \in Z$ such that $\operatorname{mult}_z(Z) = 2$ and let $k \in \{-1,...,N-2 \}$. Then, for general $D \in \mathbb{G}(k,N)$, we have:
$$ \operatorname{mult}_z P(Z,D) \leq 2.$$
\end{cor}

\begin{proof} The result is obvious for $k=-1$, since in this case $P(Z,D) = Z$. Assume that $k \geq 0$ and let $D \in \mathbb{G}(k,N)$ be a general $k$-plane. Let $u \in D$ be a general point ans let $\pi_u$ be the projection from $u$. Then, the projections
$$\pi_u|_{P(Z,u)} : P(Z,u) \rightarrow \pi_u (P(Z,u))$$ 
and
$$\pi_u|_{P(Z,D)} : P(Z,D) \rightarrow \pi_u (P(Z,D))$$
are locally isomorphisms around $z$. Moreover, we have the equality (see \cite{te}):
$$ \pi_u (P(Z,D)) = P(\pi_u (P(Z,u)), \pi_u(D)).$$
As a consequence, it is sufficient to prove the result for $k=0$. But in this case, this is an obvious application of the above lemma. Indeed, for general $u \in \mathbb{P}^N$, we have:
$$\operatorname{mult}_z P(f_Z,u) = \operatorname{mult}_z(Z)-1 = 1.$$

\end{proof}

We also need the following result.

\begin{prop} \label{multicompo} Let $X \subset \mathbb{P}^N$ be an irreducible projective variety such that $X^*$ is a hypersurface. Let $[h] \in X^*$ be such that $\operatorname{Tan}(H,X)$ has $m$ components (some of which may be embedded components), then there exists $k \in \{\-1,...,N-2 \}$ , such that for general $D \in G(k,N)$, we have:
$$ \operatorname{mult}_{[h]} P(X^*,D) \geq m.$$
\end{prop}

\begin{proof}
We only prove the result when $\operatorname{Tan}(H,X)$ is reduced and pure dimensional. The general case is done using the same ideas, it is only more technical.
\bigskip

Assume that:
$$ \operatorname{Tan}(H,X) = Y_1 \cup ... \cup Y_m,$$ where the $Y_i$ have the same codimension, say $c$. Let $D \subset {\mathbb{P}^N}^*$ be a general $\mathbb{P}^{N+1-c}$. Then, we have:
$$ \pi_D (P(X^*,D)) = (D^{\bot} \cap X)^*,$$ where $\pi_D$ is the projection from $D$. Moreover, we have $[h] \in P(X^*,D)$ and:
$$\operatorname{Tan}(D^{\bot} \cap H, D^{\bot} \cap X) = D^{\bot} \cap \operatorname{Tan}(H,X).$$
As a consequence, $ \operatorname{Tan}(D^{\bot} \cap H, D^{\bot} \cap X)$ is a $0$-dimensional scheme of degree at least $m$. In this case, it is clear that: $$\operatorname{mult}_{\pi_D([h])} \pi_D(P(X^*,D)) \geq m.$$

On the other hand, since $D$ is general, the morphism:
$$ \pi_D : P(X^*,D) \rightarrow \pi_D(P(X^*,D))$$ is locally an isomorphism around $[h]$, so that:
$$ \operatorname{mult}_{[h]} P(X^*,D) \geq m.$$

\end{proof}

Now, we can dive into the proof of \ref{multideux}.

\begin{proof} Let $T_1 \cup...\cup T_m$ be the decomposition of $\operatorname{Tan}(H,X)$ into irreducible components. If $m \geq 3$, then proposition \ref{multicompo} implies that $\operatorname{mult}_{[h]}(X^*) \geq 3$, this is impossible, so that $m \leq 2$. 

\bigskip

Assume that $m=2$. The proof of proposition \ref{multicompo} shows that these two irreducible components are scheme-theoretically linear spaces.

\bigskip

Assume that $m=1$ and let $k \in \{-1,...,N-2 \}$ such that $T_1$ is dual to some irreducible components of the reduced space underlying $\mathcal{C}_{[h]} P(X^*,D)$, for general $D \in \mathbb{G}(k,N)$. By corollary \ref{cone}, the cone $\mathcal{C}_{[h]} P(X^*,D)$ is either a hyperquadric or a linear space. Assume that it is an irreducible hyperquadric. If $k \geq 0$, we know by theorem \ref{lete} that $|\mathcal{C}_{[h]}(X^*)|^*$ is the reduced space underlying some embedded component of $\operatorname{Tan}(H,X)$. Taking $q = \dim \operatorname{Tan}(H,X)$ general hyperplane sections of $\operatorname{Tan}(H,X)$ passing through $|\mathcal{C}_{[h]}(X^*)|^*$, we see as in the proof of proposition \ref{multicompo} that for general $D' \in \mathbb{G}(q-1,N)$, we have:
$$ \operatorname{mult}_{[h]}(P(X^*,D')) \geq 3.$$ This is impossible by corollary \ref{cone}. Thus, if $\mathcal{C}_{[h]} P(X^*,D)$ is an irreducible hyperquadric, then $k=-1$, and we are in the case 1 of the proposition.

\bigskip

Finally, if $\mathcal{C}_{[h]} P(X^*,D)$ is a the union of two linear spaces or a unique linear space, then we are in case 2 or 3 of the proposition. This concludes the proof of proposition \ref{multideux}.

\end{proof}

\end{section}
\end{appendix}

\newpage

\bibliographystyle{alpha}

\bibliography{bibli}

\end{document}